\documentclass{article}
\usepackage{amssymb, amsthm, amsmath, amsfonts}
\usepackage{graphicx}
\usepackage[T2A]{fontenc}
\usepackage{float}
\usepackage{subfigure}
\usepackage{enumitem}

\newcommand{\rp}{\mathbb R \ensuremath{\text{P}}}

\newtheorem{theorem}{Theorem}
\newtheorem{lemma}{Lemma}
\newtheorem{proposition}{Proposition}

\begin{document}
\sloppy
\title{Topology of Ambient 3-Manifolds of Non-singular Flows with Twisted Saddle Orbit}
\date{}
\author{O. V. Pochinka, D. D. Shubin}

\maketitle

\begin{center}
	{\it National Research University Higher School of Economics}
\end{center}
	
%\tableofcontents
\begin{abstract}
In the present paper, non-singular Morse-Smale flows on closed orientable 3-manifolds  under the assumption that among the periodic orbits of the flow there is only one saddle one and it is twisted are considered. An exhaustive description of the topology of such manifolds is obtained. Namely, it has been established that any manifold admitting such flows is either a lens space, or a connected sum of a lens space with a projective space, or Seifert manifolds with base homeomorphic to sphere and three singular fibers. Since the latter are simple manifolds, the result obtained refutes the result that among simple manifolds, the considered flows admit only lens spaces.  
\end{abstract}
	
\section{Introduction and formulation of results}

In the present paper, we consider {\it NMS-flows} $f^t$, that is, {\it non-singular} (without fixed points) Morse-Smale flows defined on closed orientable 3-manifolds $M^3$. The nonwandering set of such flow consists of a finite number of periodic hyperbolic orbits. It is known from Azimov's work \cite{Azimov} that the ambient manifold in this case is the union of circular handles. However, in the case of a small number of periodic orbits, the topology of the manifold can be significantly refined. 
For example, only lens spaces are ambient for NMS-flows with exactly two periodic orbits.
Moreover, in \cite{PoSh} it is proved that for every lens space there are exactly two equivalence classes of such flows, except for the 3-sphere $\mathbb S^3$ and the projective space $\mathbb RP^3$, on which there is one equivalence class.

In \cite{CamposCorderoMartinezAlfaroVindel}, it is stated that the lens space is also the only {\it simple} (homeomorphic to $\mathbb S^2\times\mathbb S^1$ or {\it irreducible} -- any cylindrically embedded 2-sphere bounds the 3-ball) 3-manifold which is ambient for NMS-flows with unique saddle periodic orbit. However, this is incorrect. In the previous works of the authors \cite{Shu21} NMS-flows with exactly three periodic orbits (attractive, repelling, and saddle) are constructed on a countable set of pairwise non-homeomorphic mapping tori that are not lens spaces. Moreover, in \cite{PoSh-eqv}, necessary and sufficient conditions for the topological equivalence of such flows are obtained.

In this paper, we recognize the topology of all possible orientable 3-manifolds that admit NMS-flows with exactly one saddle periodic orbit, assuming that it is {\it twisted} (its invariant manifolds are nonorientable). 

Let us proceed to the formulation of the results.

Let $M^3$ be a connected closed orientable 3-manifold, $f^t\colon M^3\to M^3$ an NMS-flow, and $\mathcal O$ its periodic orbit. In the neighborhood of a hyperbolic periodic orbit $\mathcal O$, the flow can be simply described (up to topological equivalence). Namely, there exist linear diffeomorphism of the plane, given by matrix with positive determinant and real eigenvalues with absolute value different from one and tubular neighborhood $V_{\mathcal O}$ homeomorphic to the solid torus $\mathbb D^2\times \mathbb S^1$, in which the flow is topologically equivalent to the suspension over this diffeomorphism (see, for example, \cite{Irwin}). If both eigenvalues are greater (less than) one in absolute value, then the corresponding periodic orbit is {\it repelling (attractive)}, otherwise it is {\it saddle}. In this case, a saddle orbit is called {\it twisted} if both eigenvalues are negative and {\it untwisted} otherwise.

Let us choose {\it meridian}  $M_{\mathcal O}$  (a null-homotopic curve on $V_{\mathcal O}$ and essential on $T_{\mathcal O}$) on the torus $T_{\mathcal O}=\partial V_{\mathcal O}$  and the longitude $L_{\mathcal O}\subset T_{\mathcal O}$ (the curve homologous in the $V_{\mathcal O}$ to the orbit $\mathcal O$). We assume that the meridian $M_{\mathcal O}$ is oriented so that the pair of oriented curves $M_{\mathcal O}, L_{\mathcal O}$ defines the outer side of the solid torus boundary. Thus, the homotopy types $\langle L_{\mathcal O}\rangle=\langle1,0\rangle,\, \langle M_{\mathcal O}\rangle=\langle0,1\rangle$ of knots $L_{\mathcal O},\, M_{\mathcal O}$ are generators of the homotopy types $\langle K\rangle$ of oriented knots $K$ on the torus $T_{\mathcal O}$:
\begin{equation}\label{1}
\langle K\rangle=\langle l_{\mathcal O},m_{\mathcal O}\rangle =l_{\mathcal O}\langle L_{\mathcal O}\rangle+ m_{\mathcal O}\langle M_{\mathcal O}\rangle,
\end{equation}
where $l_{\mathcal O},\,m_{\mathcal O}\in\mathbb Z$ is the number of twists of the oriented knot $K$ around the parallel and the meridian, respectively.

\begin{figure}[H]
	\centering
	\begin{minipage}[b]{0.495\textwidth}
\center{\includegraphics[width=1\linewidth]{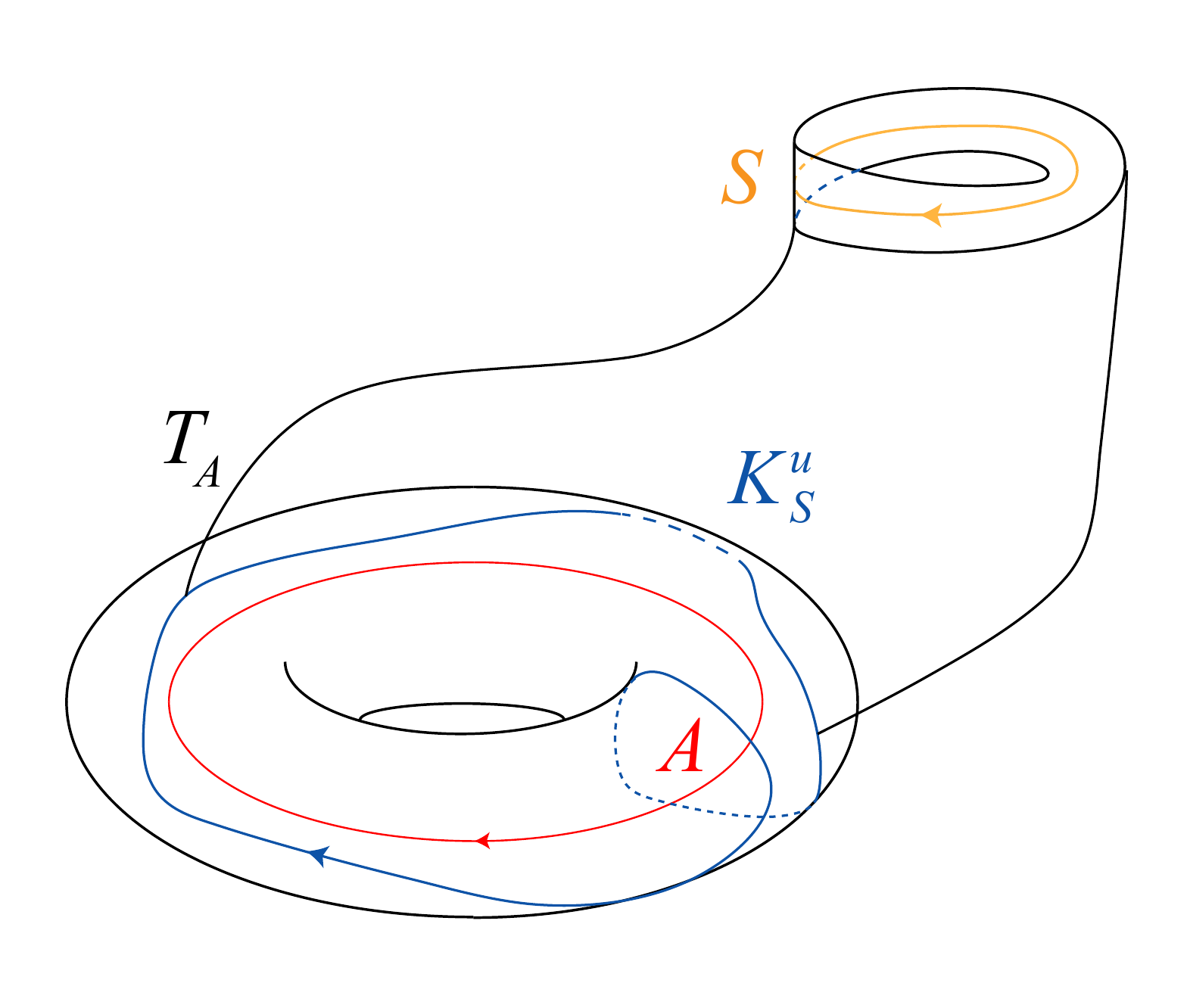}}
		\\
		{essential}
	\end{minipage}
	\hfill
	\begin{minipage}[b]{0.495\textwidth}
\center{\includegraphics[width=1\linewidth]{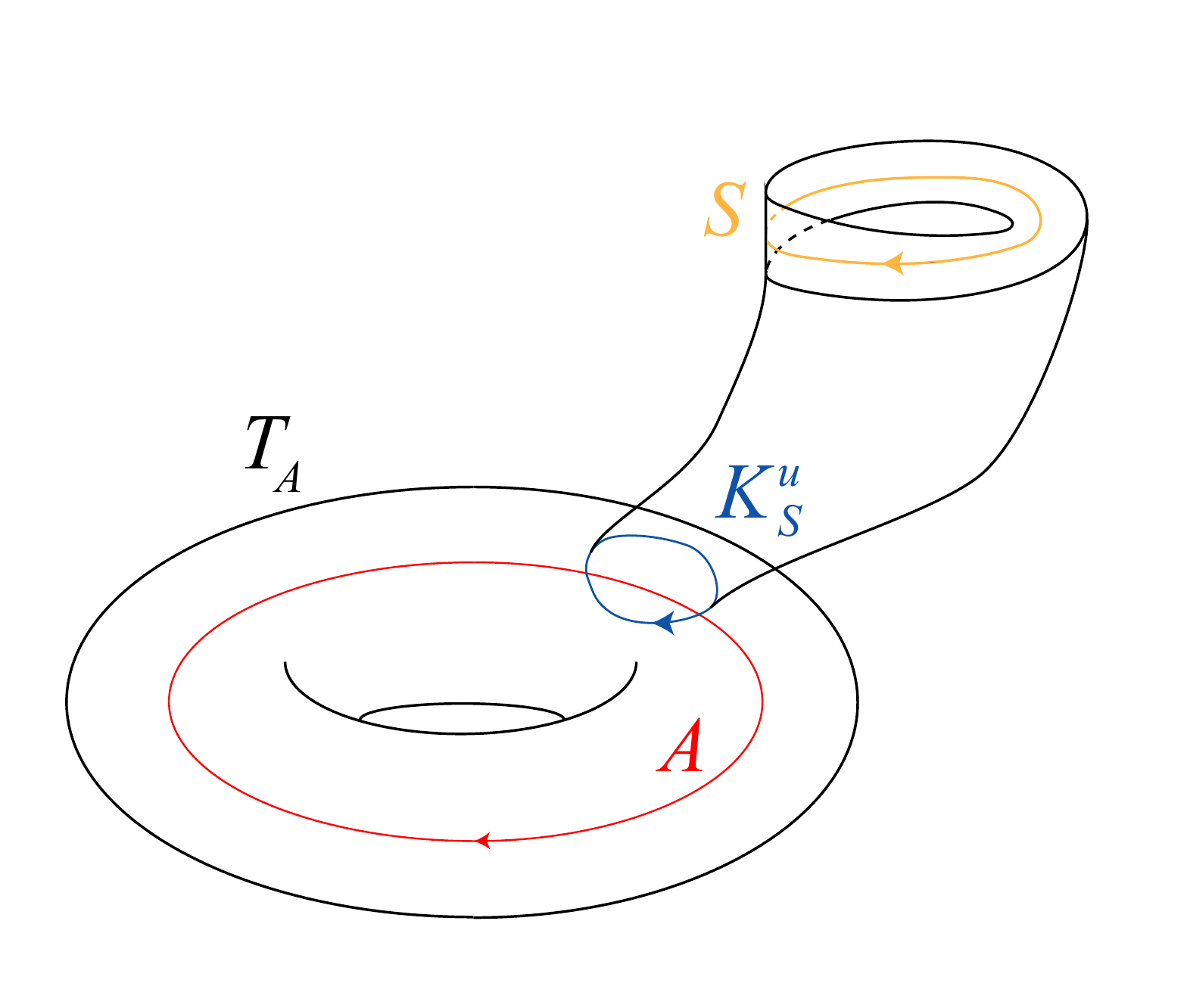}}
		\\
		{inessential}
	\end{minipage}
	\caption{Knot $K_S^u$}\label{pm-}
\end{figure}

Consider the class $G^-_3(M^3)$ of NMS-flows $f^t\colon M^3\to M^3$ with unique saddle orbit, assuming that it is twisted. Since the ambient manifold $M^3$ is the union of the stable (unstable) manifolds of all its periodic orbits, then the flow $f^t\in G^-_3(M^3)$ must have at least one attracting and at least one repelling orbit . In the Section~\ref{dyd} , we will prove the following fact.
\begin{lemma}\label{RAS} The nonwandering set of any flow $f^t\in G^-_3(M^3)$ consists of exactly three periodic orbits $S,A,R$, saddle, attracting and repelling, respectively.
\end{lemma}
Since the flow $f^t$ in the neighborhood of a periodic orbit is equivalent to a suspension over the linear diffeomorphism, the stable and unstable manifolds of these orbits have the following topology:
\begin{itemize}
\item $W^u_S\cong W^s_S\cong\mathbb R\tilde\times\mathbb S^1$ (open M\"obius strip);
\item $W^s_A\cong W^u_R\cong\mathbb R^2\times\mathbb S^1$;
\item $W^u_A\cong W^s_R\cong\mathbb S^1$.
\end{itemize} 
This fact and Lemma~\ref{RAS} immediately imply  the following proposition (for more details~see~\cite{PoSh-eqv}). 
\begin{proposition}\label{prop} 
The ambient manifold $M^3$ of any flow $f^t\in G^-_3(M^3)$ is represented as the union of three solid tori:
$$M^3 = \mathcal V_A \cup V_S \cup \mathcal V_R$$ 
with disjoint interiors being tubular neighborhoods of $A,S,R$ orbits, respectively, with the following properties:
\begin{itemize}
\item $T_S=\partial V_S$ is the union of tubular neighborhoods $T^u_S,\,T^s_S$ of knots $K^u_S=W^u_S\cap T_S,\,K^s_S=W^s_S\cap T_S $, respectively, such that $T^u_S\cap T^s_S=\partial T^u_S\cap\partial T^s_S$;
\item the torus $\mathcal T_A=\partial \mathcal V_A$ is the union of the annuli $T^u_S,\,\mathcal T$ with disjoint interiors, and the knot $K^u_S$  has homotopy type
$$\langle K^u_S\rangle=\langle l_A,m_A\rangle$$
with respect to generators $L_A,\,M_A$
\item the torus $\mathcal T_R=\partial \mathcal V_R$ is the union of the annuli $T^s_S,\,\mathcal T$ with disjoint interiors, and the knot $K^s_S$ has homotopy type
$$\langle K^s_S\rangle=\langle l_R,m_R\rangle$$
with respect to generators $L_R,\,M_R$.
\end{itemize}
\end{proposition}
Thus the knots $K^u_S\subset\mathcal T_A,\, K^s_S\subset\mathcal T_R$ are either both inessential or both essential (see~Fig.~\ref{pm-}).
Given the flow $f^t\in G^-_3(M^3)$, we define a quadruple of integers $$C_{f^t}=(l_1,m_1,l_2,m_2)$$ as follows:
\begin{itemize}
\item if the knots $K^u_S, K^s_S$ are essential on the tori $\mathcal T_A,\mathcal T_R$, then, $$C_{f^t}=(l_R,m_R,l_A,m_A);$$
\item if the knots $K^u_S, K^s_S$ are inessential on the tori $\mathcal T_A,\mathcal T_R$ %and the disc bounded by the knot $K^u_S$ on the torus $\mathcal T_A$, 
, then,
$$C_{f^t}=(0,2,l_2,m_2),$$ where $\langle l_2,m_2\rangle$ is the homotopy type of the knot on the torus $\mathcal T_R$ which is the meridian on the torus $\mathcal T_A$.
\end{itemize}
The main result of the paper is the following theorem (all the necessary information about the objects mentioned below is given in the Section~\ref{Seif}).
\begin{theorem}\label{th:top} 
	Ambient manifolds of the flows in $G^-_3(M^3)$ are lens spaces $L_{p,q}$, all connected sums of the form $L_{p,q}\# \mathbb RP^3$ and all Seifert manifolds of the form $M(\mathbb S^2, (2, 1),(\alpha_1,\beta_1),(\alpha_2,\beta_2))$. Namely, let the flow $f^t\in G^-_3(M^3)$ correspond to the collection $C_{f^t} =(l_1, m_1, l_2, m_2)$. Then
\begin{enumerate}[label={\arabic*)}]
	\item if $l_1 = 0$ and $l_2\neq 0$, then $M^3$ is homeomorphic to the manifold $L_{l_2, m_2} \#\ \rp^3$;
	\item if $l_1\neq 0$ and $l_2 = 0$, then $M^3$ is homeomorphic to the manifold $L_{l_1, m_1} \#\ \rp^3$;
	\item if $l_1 = 0$ and $l_2 = 0$, then $M^3$ is homeomorphic to $\mathbb S^2\times \mathbb S^1 \#\ \rp^3$;
	\item if $|l_1|=1$ and $|l_2|>1$, then $M^3$ is homeomorphic to the lens space $L_{p,q}$, where $p=2m_2-l_2,\,q=m_2 $;
	\item if $|l_2|=1$ and $|l_1|>1$, then $M^3$ is homeomorphic to the lens space $L_{p,q}$, where $p=2m_1-l_1,\,q=m_1 $;
	\item if $|l_1l_2|=1$, then $M^3$ is homeomorphic to the sphere $\mathbb S^3$;
	\item if $|l_1|>1$ and $|l_2|>1$ then $M^3$ is homeomorphic to the simple Seifert manifold $M(\mathbb S^2,(2,1), (l_1,\beta_1) ,(l_2,\beta_2)),\,\beta_im_i\equiv 1 \pmod{l_i},\ i=1,2,$ and is not homeomorphic to any lens space.
\end{enumerate}
\end{theorem}

{\it {\bf Acknowledgments:} The author is partially supported by Laboratory of Dynamical Systems and Applications NRU HSE, grant of the Ministry of science and higher education of the RF, ag. \textnumero\ 075-15-2022-1101

}

\section{Necessary information on the topology of 3-manifolds}\label{Seif}
\subsection{Lens spaces}
Everywhere below, we assume that generators of homotopy types of knots on the boundary $\partial\mathbb V$ of the standard solid torus $\mathbb V = \mathbb D^2\times \mathbb S^1$ are the meridian $\mathbb M=(\partial \mathbb D^2) \times \{y\},\ y\in\mathbb S^1$ with homotopy type $\langle 0,1\rangle$ and parallel $\mathbb L=\{x\}\times \mathbb S^1,\,x\in \partial \mathbb D^2$ with homotopy type $\langle 1,0\rangle$.

\emph{Lens space} is a three-dimensional manifold $L_{p,q} = V_1 \cup_j V_2$, which is the result of gluing together two copies of the solid torus $V_1=\mathbb V,\, V_2=\mathbb V$ by some homeomorphism $j \colon \partial V_1\to \partial V_2$ such that $j_*(\langle 0,1\rangle)=\langle p,q\rangle$.
\begin{proposition}[\cite{Rolfsen}]\label{lens-class}
	Two lens spaces $L_{p,q},\,L_{p',q'}$
	are homeomorphic (up to preserving the numbering of copies) if and only if $p = \pm p',\ q \equiv \pm q' \pmod{|p|}$.
\end{proposition}
\subsection{Dehn surgery along knots and links}
Suppose the following data are given:
\begin{enumerate}[label=(\alph*)]
	\item a closed 3-manifold $M$;
	\item knot $\gamma\subset M$;
	\item tubular neighborhood $U_\gamma$ of $\gamma$ with standard generators on $\partial U_\gamma$: meridian $M_\gamma$ and longitude $L_\gamma$;
	\item homeomorphism $h\colon\partial\mathbb V\to\partial U_\gamma$ inducing an isomorphism defined by the matrix 
	$h_{*}=\begin{pmatrix}
		\xi & -\nu \\
		\beta & \alpha
	\end{pmatrix}\in SL(2,\mathbb Z)$ with respect to given generators.
\end{enumerate}
Manifold
$$M_{\gamma,h} = (M \setminus \mathrm{int }\ U_\gamma) \cup_{h} \mathbb V$$ is called the \emph{ manifold obtained from $M$ by Dehn surgery along the knot $\gamma$}.
\begin{proposition}[\cite{Rolfsen}]\label{K}
	If there exists a homeomorphism $\varphi\colon M\to M$ such that $\varphi(\gamma)=\gamma'$ and the matrices $h_{*}=\begin{pmatrix}
		\xi & -\nu \\
		\beta & \alpha
	\end{pmatrix},\,h'_{*}=\begin{pmatrix}
		\xi' & -\nu' \\
		\beta' & \alpha'
	\end{pmatrix}$ are related by $\beta=\beta',\,\alpha\equiv\alpha'\pmod {\beta}$, then $M_{\gamma,h}\cong M_{\gamma', h'}$.
\end{proposition}
Thus, up to homeomorphism, the topology of the manifold $M_{\gamma,h}$ is determined by the equivalence class of the knot $\gamma$ and a pair of coprime numbers $\beta,\alpha$ called the {\it surgery coefficients of the knot $\gamma$} . So, the manifold
$M_{\gamma,h}$ will be written further as $M_{\gamma}$, implying that $\gamma$ is an equipped knot.

Naturally, the manifold $M$ is restored from $M_{\gamma,h}$ by inverse surgery. Namely, we denote by $p_{\gamma,h}:(M \setminus \mathrm{int }\ U_\gamma) \sqcup \mathbb V\to M_{\gamma,h}$ the natural projection.
Let $\tilde\gamma=p_{\gamma,h}(\{0\}\times\mathbb S^1),\,U_{\tilde\gamma}=p_{\gamma,h}(\mathbb V ),\,\tilde h=p_{\gamma,h}h^{-1}:\partial U_\gamma\to\partial U_{\tilde\gamma}$. Then
\begin{equation}\label{obr}
	M=(M_{\gamma,h})_{\tilde\gamma,\tilde h}.
\end{equation}
The following assertions follow directly from the relation (\ref{obr}) and the Proposition~\ref{K}. 
\begin{proposition}\label{ts} Let $\gamma\subset M$ be equipped with $\beta,\alpha$. Then $$M=(M_{\gamma})_{\tilde\gamma},$$ where $\tilde\gamma$ is equipped with $-\beta,\xi$ satisfying $\xi\alpha+\nu\beta=$1.
\end{proposition}
\begin{proposition}\label{lm} Let $L_{p,q} = V_1 \cup_j V_2$, where $j_*(\langle 0,1\rangle)=\langle p,q\rangle$ and $\mathbb S^3=V_1 \cup_{j_0} V_2$, where $j_{0*}(\langle 0,1\rangle)=\langle 1,0\rangle$. Then $$L_{p,q}\cong\mathbb S^3_{M_1},$$ where $M_1$ is the meridian of the torus $V_1$ with framing $q,p$.
\end{proposition}
Dehn surgery are naturally generalized to the case where $\gamma = \gamma_1\sqcup\dots\sqcup \gamma_r\subset M$ is a disjoint union (link) of equipped knots.
The resulting manifold $M_\gamma$ in this case is called the \emph{ manifold obtained from the manifold $M^3$ by Dehn surgery along the equipped link $\gamma$}.
A link $\gamma=\gamma_1\sqcup\dots\sqcup \gamma_r\subset M$ is called {\it trivial} if knots $\gamma_1,\dots,\gamma_r$ bound pairwise disjoint 2-discs $d_1,\dots ,d_r\subset M$.
\begin{proposition}[\cite{Rolfsen}]\label{prop:dehn-connected-sum}
	Let $\gamma = \gamma_1\sqcup\dots\sqcup \gamma_r\subset M$ be a trivial link equipped with  $q_1,p_1,\dots,q_r,p_r$. Then
	$$M_{\gamma} \cong M\#L_{p_1,q_1}\# \dots L_{p_r,q_r}.$$
\end{proposition}

\subsection{Seifert fiber space} 
A solid torus $\mathbb V$ split into fibers of the form $\{x\}\times \mathbb S^1$ is called a \emph{trivially fibred solid torus}. Consider the solid torus $\mathbb V = \mathbb D^2\times \mathbb S^1$ as the cylinder $\mathbb D^2\times [0, 1]$ with the bases glued due to the $2\pi\nu/ \alpha$ angle rotation for coprime integers $\alpha, \nu,\ \alpha > 1$. The partition of the cylinder into segments of the form $\{x\}\times [0, 1]$ determines the partition of this solid torus into circles called \emph{fibers}. The segment $\{0\}\times [0, 1]$ generates a fiber which we call {\it exceptional}, all other ({\it ordinary}) fibers of the solid torus wrap $\alpha$ times around the exceptional fiber and $\nu $ times around the solid torus meridian. The number $\alpha$ is called the \emph{multiplicity} of the exceptional fiber.
A solid torus with such a partition into fibers is called a \emph{nontrivially fibred solid torus} with {\it orbital invariants} $(\alpha,\nu)$.

\emph{Seifert manifold} is a compact, orientable 3-manifold $M$ split into disjoint simple closed curves (fibers) in such a way that each fiber has a neighborhood consisting of layers, fiberwise homeomorphic to a fibred solid torus. Such a partition is called \emph{Seifert fibration}. The fibers which under such homeomorphism correspond to the center of a non-trivially fibred solid torus are called \emph{exceptional}.

Two Seifert fiberings $M,M'$ are called \emph{isomorphic} if there exists a homeomorphism $h\colon M\to M'$ such that the image of each fiber of one bundle is a fiber of the second bundle. It is easy to show (see, for example, \cite[Proposition 10.1]{MatFom}) that two bundles of a solid torus with orbital invariants $(\alpha_1,\nu_1);\,(\alpha_2,\nu_2)$ are isomorphic  if and only if $\alpha_1=\alpha_2\,(=\alpha);\,\nu_1 \equiv \nu_2 \pmod{\alpha}$.

The \emph{base} of a Seifert manifold $M$ is a compact surface $\Sigma = M/_\sim$, where $\sim$ is an equivalence relation such that $x \sim y$ if and only if $x$ and $y$ belong to the same fiber. It is easy to show (see, for example, \cite[Proposition 10.2]{MatFom}) that the base of any solid torus bundle is a disc. The base of any Seifert manifold is a compact surface, and Seifert bundles with non-homeomorphic bases are not isomorphic (see, for example, \cite{MatFom}).

Thus, any Seifert fibering $M$ with a given base $\Sigma$ and orbital invariants
$(\alpha_1,\nu_1),\dots,(\alpha_r,\nu_r),\,r\in\mathbb N$
obtained from the manifold $\Sigma\times\mathbb S^1$ by Dehn surgery along the link $\gamma=\bigsqcup\limits_{i=1}^r \gamma_i$, where $\gamma_i=\{s_i\}\times\mathbb S^1,\,s_i\in\Sigma$ is a knot with coefficients $\beta_i,\alpha_i,\,\nu_i\beta_i\equiv 1\pmod{\alpha_i}$. Therefore, the conventional notation for such a Seifert fibration is
$$M(\Sigma, (\alpha_1, \beta_1),\dots, (\alpha_r, \beta_r)).$$
\begin{proposition}[\cite{Hatcher}, \cite{MatFom}]\label{Seifert-class}
	Seifert fibrations $M(\Sigma, (\alpha_1, \beta_1),\dots, (\alpha_r, \beta_r))$ and $M'(\Sigma', (\alpha'_1, \beta'_1),\ dots, (\alpha'_{r'}, \beta'_{r'}))$ are isomorphic (preserving orientation of fibers) if and only if the following conditions are satisfied:
	\begin{itemize}
		\item $\Sigma$ is homeomorphic to $\Sigma'$;
		\item $r = r'$; $\alpha_i = \alpha'_i;\,\beta_i \equiv \pm\beta'_i \pmod{\alpha_i}$ for $i\in\{1,\dots, r\}$;
		\item if the surface $\Sigma$ is closed, then $\sum\limits_{i=1}^{r} \frac{\beta_i}{\alpha_i} = \sum\limits_{i=1}^{r}\frac{\beta'_i}{\alpha'_i}$.
	\end{itemize}
\end{proposition}
\begin{proposition}[\cite{Hatcher}, Proposition 1.12]\label{prosto} All closed orientable Seifert manifolds are prime except $M(\mathbb S^2,(2,1),(2,1) ,(2,1),(2,1))\cong\mathbb RP^3\#\mathbb RP^3$.
\end{proposition}
\begin{proposition}[\cite{GaigesLange}]\label{lens-Seifert} A 3-manifold admits a Seifert fibration with base sphere and at most two singular fibers if and only if it is homeomorphic to lens space. Wherein,
	\begin{itemize}
		\item the only manifold which admits fibering without singular fibers admits is $\mathbb S^2\times\mathbb S^1$;
		\item $M(\mathbb S^2, (\alpha, \beta))\cong L_{\beta,\alpha}$;
		\item $M(\mathbb S^2, (\alpha_1, \beta_1), (\alpha_2, \beta_2))\cong L_{p,q}$, where $p=\beta_1\alpha_2-\alpha_1\beta_2 ,\,q=\beta_1\nu_2+\alpha_1\xi_2$ and $\alpha_2\xi_2+\nu_2\beta_2 = 1$.
	\end{itemize}
\end{proposition}
It follows from the above statement, in particular, that any lens space admits more than a unique Seifert fibration. However, as the result below shows, any such fibration with base sphere cannot have more than two singular fibers.
\begin{proposition}[\cite{MatFom}]\label{3net}
No lens space admits a Seifert fibration with base sphere and more than two singular fibers.
\end{proposition}

\section{Dynamicrs of the flows $G^-_3(M^3)$}\label{dyd}This section is devoted to the proof of the Lemma~\ref{RAS}: the nonwandering set of any flow $f^t\in G^-_3(M^3)$ consists of exactly three periodic orbits $S,A,R$, saddle, attracting and repelling, respectively.
\begin{proof} The basis of the proof is the following representation of the ambient manifold $M^3$ of the NMS-flow $f^t$ with the set of periodic orbits $Per_{f^t}$ (see, for example, \cite{Sm})
\begin{equation}\label{Mob}
	M^3 = \bigcup\limits_{\mathcal O \in Per_{f^t}} W^u_{\mathcal O}=\bigcup\limits_{\mathcal O \in Per_{f^t}} W^s_{\mathcal O},
\end{equation} as well as the asymptotic behavior of invariant manifolds
\begin{equation*}\label{neust}
	cl(W^u_{\mathcal O}) \setminus W^u_{\mathcal O} = \bigcup\limits_{\tilde{\mathcal O} \in Per_{f^t}\colon W^u_{\mathcal O}\cap W^s_{\mathcal O}\neq \varnothing} W^u_{\tilde{\mathcal O}},
\end{equation*}
\begin{equation*}\label{ust}
	cl(W^s_{\mathcal O}) \setminus W^s_{\mathcal O} = \bigcup\limits_{\tilde{\mathcal O} \in Per_{f^t}\colon W^s_{\mathcal O}\cap W^u_{\mathcal O}\neq \varnothing} W^s_{\tilde{\mathcal O}}.
\end{equation*}
In particular, it follows from the above relations that any NMS-flow has at least one attracting orbit and at least one repelling one. Moreover, if an NMS-flow has a saddle periodic orbit, then the basin of any attracting orbit has a non-empty intersection with an unstable manifold of at least one saddle orbit (see~Proposition~2.1.3~\cite{begin}) and a similar situation with the basin of a repelling orbits.

Now let $f^t\in G_3^-(M^3)$ and $S$ be its only saddle orbit. It follows from the relation~(\ref{Mob}) that $W^u_S\setminus S$ intersects only basins of attracting orbits. Since the set $W^u_S\setminus S$ is connected and the basins of attracting orbits are open, then $W^u_S$ intersects exactly one such basin. Denote by $A$ the corresponding attracting orbit. Since there is only one saddle orbit, there is only one attracting orbit. Similar reasoning for $W^s_S$ leads to the existence of a unique repulsive orbit $R$.
\end{proof}

\section{Topology of  ambient manifolds of flows of the class $G^-_3(M^3)$}
In this section, we prove the Theorem~\ref{th:top}: flows of class $G^-_3(M^3)$ admit all lens spaces $L_{p,q}$, all connected sums of the form $L_{p, q}\#\mathbb RP^3$ and all Seifert manifolds of the form $M(\mathbb S^2, (2, 1),(\alpha_1,\beta_1),(\alpha_2,\beta_2))$. Namely, let the flow $f^t\in G^-_3(M^3)$ have the invariant $C_{f^t} =(l_1, m_1, l_2, m_2)$. Then
\begin{enumerate}[label={\arabic*)}]
	\item if $l_1 = 0$ and $l_2\neq 0$, then $M^3$ is homeomorphic to $L_{l_2, m_2} \# \rp^3$;
	\item if $l_1\neq 0$ and $l_2 = 0$, then $M^3$ is homeomorphic to $L_{l_1, m_1} \# \rp^3$;
	\item if $l_1 = 0$ and $l_2 = 0$, then $M^3$ is homeomorphic to $\mathbb S^2\times \mathbb S^1 \#\ \rp^3$;
	\item if $|l_1|=1$ and $|l_2|>1$, then $M^3$ is homeomorphic to the lens space $L_{p,q}$, where $p=2m_2-l_2,\,q=m_2 $;
	\item if $|l_2|=1$ and $|l_1|>1$, then $M^3$ is homeomorphic to the lens space $L_{p,q}$, where $p=2m_1-l_1,\,q=m_1 $;
	\item if $|l_1l_2|=1$, then $M^3$ is homeomorphic to the sphere $\mathbb S^3$;
	\item if $|l_1|>1$ and $|l_2|>1$ then $M^3$ is homeomorphic to the simple Seifert manifold $M(\mathbb S^2,(2,1), (l_1,\beta_1) ,(l_2,\beta_2)),\,\beta_im_i\equiv 1 \pmod{l_i},\ i=1,2,$ and is not homeomorphic to any lens space.
\end{enumerate}

\begin{proof} The idea of the proof is to recognize that  the sphere $\mathbb S^3$ is obtained by Dehn surgery along a link consisting of a saddle orbit $S$ and some node $\gamma$ from the carrier manifold $M^3$ of the flow $f^t$ . Then, due to the relation (\ref{obr}), $M^3\cong\mathbb S^3_{\tilde S\sqcup\tilde\gamma}$, which allows us to describe the topology of the manifold $M^3$ using the set $ C_{f^t} = (l_1, m_1, l_2, m_2)$. Let's break the discussion down into steps.

{\bf 1. Dehn surgery along a saddle orbit $S$.} Let us show that the following relation is true for a saddle orbit $S$ with coefficients $-1,0$:
$$M^3_S \cong L_{p, q}.$$
Let us put
$$\mathbb V_+=\{(d_1, d_2, s)\in\mathbb V|\ d_1 \geqslant 0\},\,\mathbb T_+ = \{(d_1, d_2, s)\in\partial \mathbb V|\ d_1 \geqslant 0\},$$
$$\mathbb V_-=\{(d_1, d_2, s)\in\mathbb V|\ d_1 \leqslant 0\},\,\mathbb T_- = \{(d_1, d_2, s)\in\partial \mathbb V|\ d_1 \leqslant 0\}.$$
Let $h\colon \partial \mathbb V\to \partial V_S$ be a homeomorphism such that
$$h(\mathbb T_+) = T^u_S,\ h(\mathbb T_-) = T^s_S.$$ Then $h_*(\langle1, 0\rangle) = \langle 2, 1\rangle$ , which means that
the matrix induced by it has the form:
$$h_* = \begin{pmatrix}
	2&1\\
	\beta & \alpha
\end{pmatrix}.$$
From the condition $2\alpha - \beta = 1$, we get that $\beta = -1,\ \alpha = 0$. Consider the Dehn surgery $M^3_S$ on $M^3$ along the knot $S$ with neighborhood $V_S$ and coefficients $-1,0$. Let $v_S\colon (M^3 \setminus \mathrm{int }\ V_S) \sqcup \mathbb V\to M^3_S$ be the natural projection. For simplicity, we keep the notation of all objects on $v_S(M^3\setminus int\,V_S)$ the same as they were on $M^3\setminus int\,V_S$ and set $\tilde S=v_S(\{ 0\}\times\mathbb S^1),\,V_{\tilde S} = v_S(\mathbb V)$. Then $M^3_S$ is the union of two filled tori $\tilde{\mathcal V}_A=\mathcal V_A\cup v_S(\mathbb V_+)$ and $\tilde{\mathcal V}_R=\mathcal V_R\cup v_S(\mathbb V_-)$ such that $\tilde{\mathcal V}_A\cap\tilde{\mathcal V}_R=\partial\tilde{\mathcal V}_A\cap\partial\tilde{\mathcal V}_R$ and hence $M^3_S \cong L_{a, b}$ for some coprime integers $a, b$.

{\bf 2. Reverse Dehn surgery on lens $L_{a, b}$ along the knot $\tilde S$.}
Let $\tilde{\mathcal T}_A=\partial\tilde{\mathcal V}_A,\,\tilde{\mathcal T}_R=\partial\tilde{\mathcal V}_R$ and $L_{a, b}=\tilde{\mathcal V}_A\cup\tilde{\mathcal V}_R$.
From the Proposition~\ref{ts} we get that $M^3=(L_{a, b})_{\tilde S},$ where $\tilde S$ is a knot with coefficients $1,2$.
For knots $\delta\subset\tilde{\mathcal T}_A(=\tilde{\mathcal T}_R)$ denote by $\langle\delta\rangle_A,\,\langle\delta\rangle_R$ the homotopy types of the knot $ \delta$ on the tori $\tilde{\mathcal T}_A,\,\tilde{\mathcal T}_R$, respectively.
Then for cases 1)-3) from the definition of $C_{f^t}$ we have the following relations.
\begin{enumerate}[label={\arabic*)}]
	\item If $l_1 = 0$ and $l_2\neq 0$, then either $\langle\tilde S\rangle_A=\langle0,0\rangle$ or $\langle\tilde S\rangle_A=\langle0,1\rangle$. In the first case $\langle\tilde S\rangle_R=\langle0,0\rangle$ and $\langle M_A\rangle_R=\langle l_2, m_2\rangle$, which means $a = l_2,\ b = m_2$. Then by the Statement~\ref{lm}, $\mathbb S^3_{M_A}=L_{l_2,m_2}$, where $M_A$ is the meridian of the torus $\tilde{\mathcal T}_A$ equipped with $m_2 ,l_2$. Thus, $M^3=\mathbb S^3_{\tilde S\sqcup M_A}$. Since the knots $\tilde S\sqcup M_A$ form a trivial link on the sphere $\mathbb S^3$ ($M_A$ can be chosen not to intersect $\tilde S$ by \ref{K}), then, by virtue of Propositions~\ref{K},~and~\ref{prop:dehn-connected-sum}, $$M^3 = \mathbb S^3_{\tilde S\sqcup M_A}\cong L_{2,1}\#L_{l_2,m_2}=L_ {l_2,m_2}\#\rp^3.$$
	Similarly, if $\langle\tilde S\rangle_A=\langle0,1\rangle$ then $\langle M_A\rangle_R = \langle\tilde S\rangle_R=\langle l_2, m_2\rangle$, and so $a = l_2,\ b = m_2$. Since $M_A$ can also be chosen to be disjoint from $\tilde S$, then $$M^3\cong L_{l_2,m_2}\#\rp^3.$$
	\item If $l_1 \neq 0$ and $l_2 = 0$, then $\langle \tilde S\rangle_R = \langle 0, 1\rangle$. Then $\langle M_R\rangle_A = \langle\tilde S\rangle_A=\langle l_1, m_1\rangle$, and hence $a = l_1,\ b = m_1$, whence, from arguments similar to the above, we obtain
	$$M^3\cong L_{l_1,m_1}\#\rp^3.$$
	\item If $l_1 =l_2 = 0$, then $\langle \tilde S\rangle_R = \langle 0, 1\rangle$. Then $\langle M_R\rangle_A = \langle\tilde S\rangle_A=\langle 0, 1\rangle$, and hence $a = 0,\ b = 1$, whence it follows that
	$$M^3\cong \rp^3 \# L_{0, 1}= \rp^3 \# \mathbb S^2\times\mathbb S^1.$$
\end{enumerate}

{\bf 3. Seifert fibration on the manifold $M^3$.} To prove the remaining points, we note that in the case when $l_1l_2\neq 0$, the manifold $M^3=\mathcal V_A\cup V_S\cup\mathcal V_R$ has a Seifert fibration. Indeed, in this case the fibration $V_S$ of the solid torus with singular fiber $S$ and orbital invariants $(2,1)$ contains the knots $K^u_S,\,K^s_S$ as fibers and extends to the solid torus $\mathcal V_A, \ \mathcal V_R$ fibration with fibers $A,\,R$ (possibly not singular) and orbital invariants $(l_1, m_1)$, $(l_2, m_2)$, respectively. In this way,
\begin{equation*}\label{eq:Seifert-proof-0}
	M^3_S\cong M(\Sigma, (l_1, \beta_1), (l_2, \beta_2))\cong L_{a, b},\quad \beta_im_i\equiv 1 \pmod{l_i}.
\end{equation*}
Let us show that the base $\Sigma$ of such a bundle is a 2-sphere.

Let $\sim$ be an equivalence relation whose equivalence classes are the fibers of this fibration. The figure \ref{base1} shows the meridian disks $D_A$, $D_S$, $D_R$ of the tori $\mathcal V_A,V_S,\mathcal V_R$, respectively, the segments containing equivalent points are marked with the same color.
\begin{figure}[!h]
\center{\includegraphics[width=0.9\textwidth]{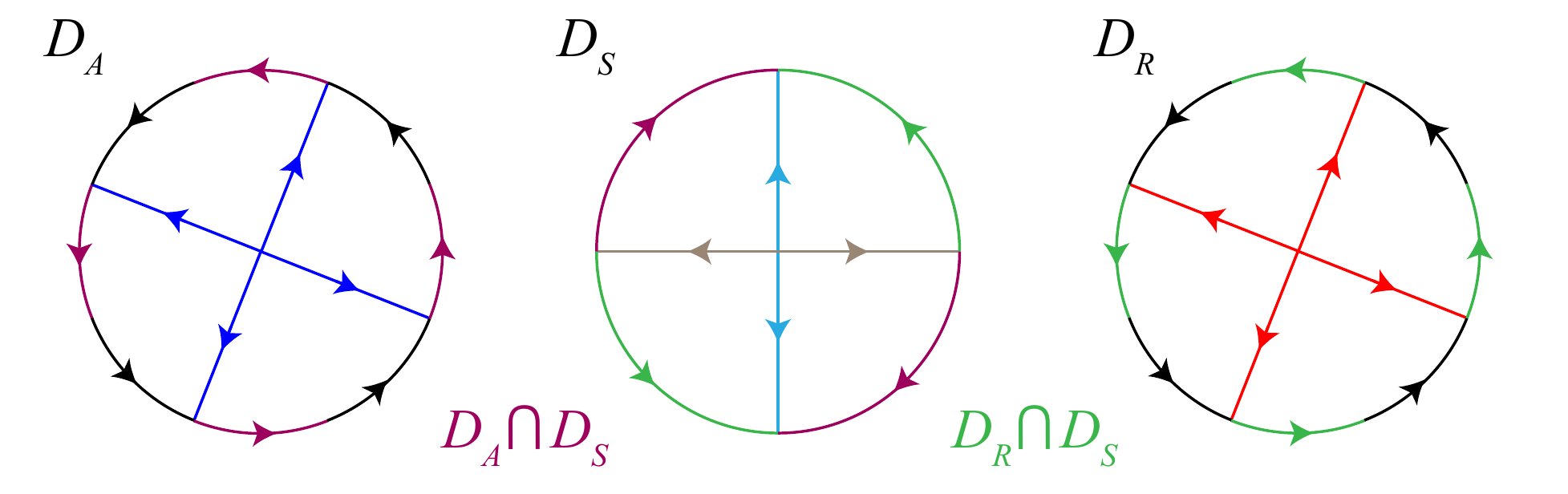}}
	\caption{Disks $D_A,D_S,D_R$}\label{base1}
\end{figure}
Gluing the equivalent points in the disks $D_A$, $D_S$, $D_R$, respectively, we obtain the disks $\hat D_A = \mathcal V_A/\sim,\ \hat D_S = V_S,\ \hat D_R = \mathcal V_R$, in which each fiber, except for the boundary layers, is represented by one point, and each boundary layer is represented by two points on different disks (see Fig.~\ref{base2}).
\begin{figure}[!h]
\center{\includegraphics[width=0.9\textwidth]{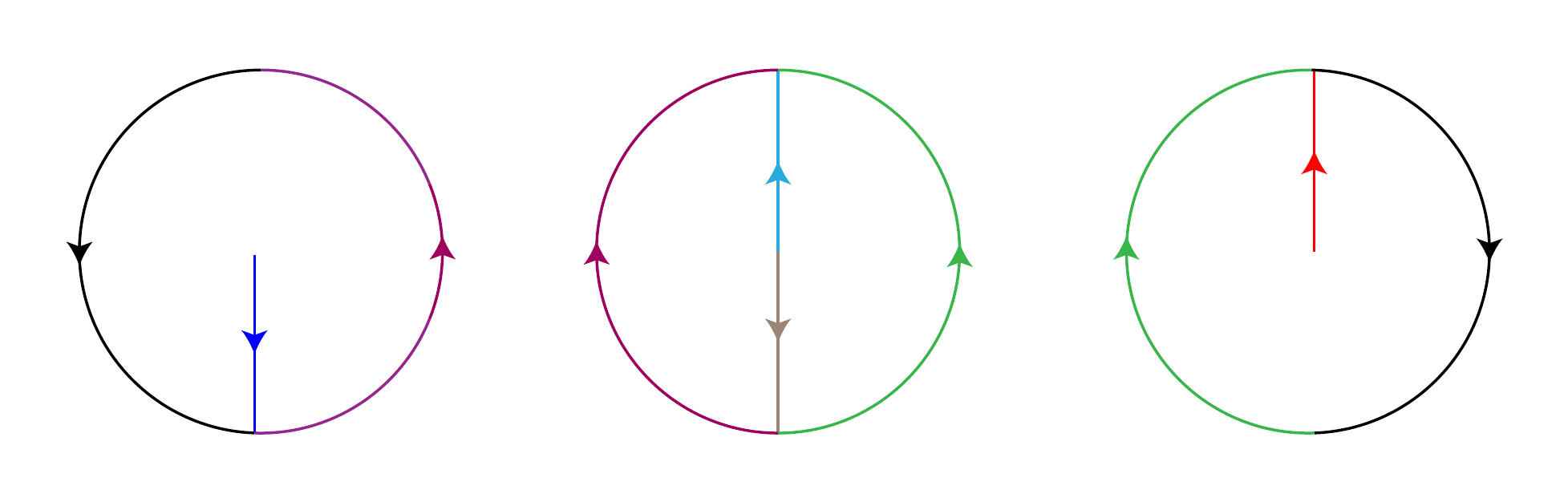}}
	\caption{Disks $\hat D_A$, $\hat D_S$, $\hat D_R$}\label{base2}
\end{figure} 
By gluing the equivalent points in the disks $\hat D_A$, $\hat D_S$, $\hat D_R$ we obtain the sphere $\mathbb S^2$ (see Fig.~\ref{base3}), which is the base of the fibration given on $M^3$.
\begin{figure}[!h]
\center{\includegraphics[width=0.9\textwidth]{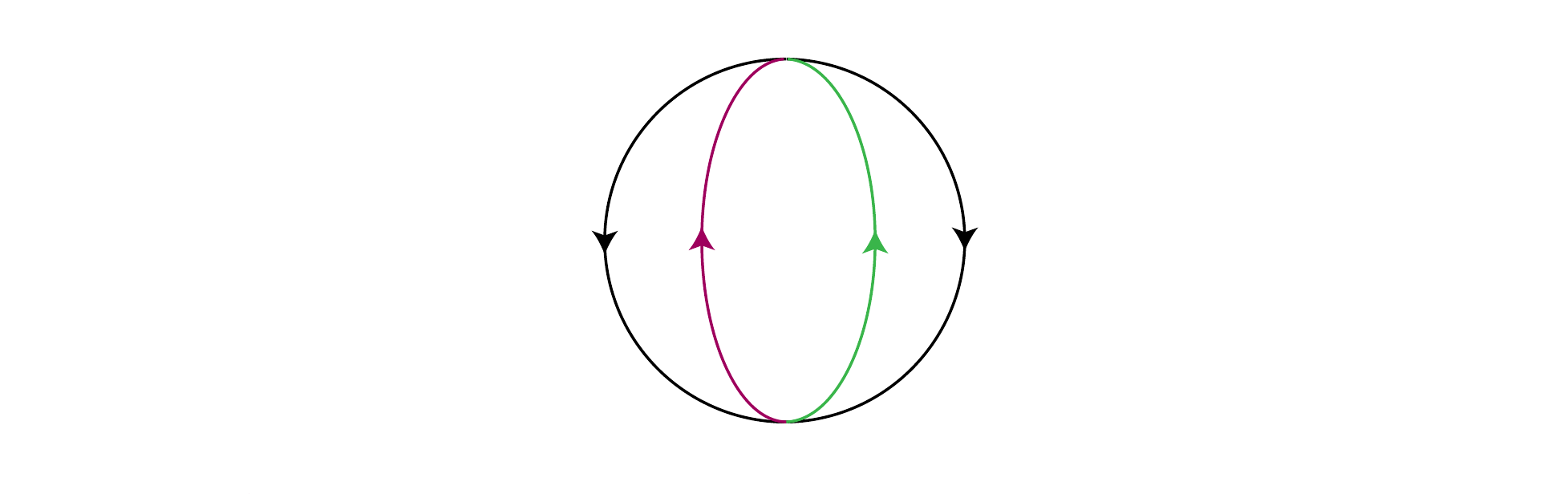}}
	\caption{$\Sigma\cong\mathbb S^2$}\label{base3}
\end{figure}
So,
\begin{equation}\label{eq:Seifert-proof}
M^3 \cong M(\mathbb S^2, (2, 1), (l_1, \beta_1), (l_2, \beta_2)),\quad \beta_im_i\equiv 1 \pmod{l_i}.
\end{equation}
\begin{enumerate}[label={\arabic*)}, resume]
	\item If $|l_1| = 1,\ |l_2| > 1$, then the fiber $A$ is ordinary and, by the Proposition~\ref{lens-Seifert}, $$M^3 \cong M(\mathbb S^2, (2, 1), (l_2, \beta_2)) \cong L_{p, q},$$
	where $p = 2m_2 - l_2,\ q = m_2$.
	\item If $|l_1| > 1,\ |l_2| = 1$, then the fiber $R$ is ordinary and, according to Proposition~\ref{lens-Seifert},
	$$M^3 \cong L_{p, q},$$
	where $p = 2m_1 - l_1,\ q = m_1$.
	\item If $|l_1l_2| = 1$, then both fibers $A,R$ are ordinary, and hence, by the Proposition~\ref{lens-Seifert}, and Proposition \ref{lens-class} $$M^3\cong M(\mathbb S^2, (2, 1))\cong L_{1, 2}\cong L_{1,0}\cong \mathbb S^3.$$
	\item If $|l_1| > 1,\ |l_2| > 1$, then $M^3$ is a Seifert manifold with three exceptional fibers
	$$M^3 \cong M(\mathbb S^2, (2, 1), (l_1, \beta_1), (l_2, \beta_2)),\quad \beta_im_i\equiv 1 \pmod{l_i}.$$ By Proposition \ref{prosto}, $M^3$ is simple, and by Proposition \ref{3net} it is not homeomorphic to lens space.
\end{enumerate}
\end{proof}

\bibliographystyle{ieeetr}
\bibliography{refs_EN}

\end{document}